\input amstex
\documentstyle{amsppt}
%%\loadbold
\magnification=\magstep1
\pageheight{9.0truein}
\pagewidth{6.5truein}
\NoBlackBoxes

\def\ZZ{{\Bbb Z}}

\def\CC{{\Bbb C}}
\def\Fract{\operatorname{Fract}}

\def\spec{\operatorname{spec}}

\def\delspec{\Delta\operatorname{-spec}}
\def\hdelspec{(H,\Delta)\operatorname{-spec}}

\def\delprim{\Delta\operatorname{-prim}}
\def\Pspec{\operatorname{P{.}spec}}
\def\Pprim{\operatorname{P{.}prim}}
\def\Pcore{\operatorname{P{.}core}}

\def\Hsymp{H\operatorname{-symp}}
\def\C{{\Cal C}}
\def\E{{\Cal E}}
\def\L{{\Cal L}}
\def\M{{\Cal M}}

\def\O{{\Cal O}}

\def\Hhat{\widehat{H}}
\def\Jhat{\widehat{J}}
\def\Der{\operatorname{Der}}
\def\kx{k^\times}
\def\kbar{\overline{k}}

\def\gfrak{{\frak g}}
\def\mfrak{{\frak m}}

\def\BrGd{{\bf 1}}
\def\BGY{{\bf 2}}
\def\BrGr{{\bf 3}}
\def\CaWe{{\bf 4}}
\def\ChOh{{\bf 5}}
\def\Dix{{\bf 6}}
\def\Dixbook{{\bf 7}}
\def\EtSc{{\bf 8}}
\def\Far{{\bf 9}}
\def\FaLe{{\bf 10}}
\def\Ful{{\bf 11}}
\def\specstrat{{\bf 12}}
\def\qaffquo{{\bf 13}}
\def\GoYa{{\bf 14}}
\def\HoLeA{{\bf 15}}
\def\HoLeB{{\bf 16}}
\def\Hum{{\bf 17}}
\def\Jos{{\bf 18}}
\def\KoSo{{\bf 19}}
\def\Moe{{\bf 20}}
\def\Oh{{\bf 21}}
\def\OhPa{{\bf 22}}
\def\OPS{{\bf 23}}
\def\Vandef{{\bf 24}}
\def\Van{{\bf 25}}

\topmatter

\title A Dixmier-Moeglin equivalence for Poisson algebras with torus actions
\endtitle

\rightheadtext{Dixmier-Moeglin equivalence for Poisson algebras}

\author K. R. Goodearl \endauthor

\address Department of Mathematics, University of California,
Santa Barbara, CA 93106, USA\endaddress
\email goodearl\@math.ucsb.edu\endemail


\thanks 
This research was partially supported by 
National Science Foundation grants DMS-9622876 and DMS-0401558. \endthanks

\abstract A Poisson analog of the Dixmier-Moeglin equivalence is established
for any affine Poisson algebra $R$ on which an algebraic torus $H$ acts
rationally, by Poisson automorphisms, such that $R$ has only finitely many prime
Poisson $H$-stable ideals. In this setting, an additional characterization of the
Poisson primitive ideals of $R$ is obtained -- they are precisely the prime
Poisson ideals maximal in their $H$-strata (where two prime Poisson ideals are in
the same $H$-stratum if the intersections of their $H$-orbits coincide). Further,
the Zariski topology on the space of Poisson primitive ideals of $R$ agrees with
the quotient topology induced by the natural surjection from the maximal ideal
space of $R$ onto the Poisson primitive ideal space. These theorems apply to many
Poisson algebras arising from quantum groups.

The full structure of a Poisson algebra is not necessary for the results of this
paper, which are developed in the setting of a commutative algebra equipped with
a set of derivations.
\endabstract

\endtopmatter

\document

\head Introduction \endhead

Motivated by existing and conjectured roles of Poisson structures in the theory
of quantum groups, we address some problems in the ideal theory of Poisson
algebras. Recall, for example, that Hodges and Levasseur \cite{\HoLeA, \HoLeB}
and Joseph \cite{\Jos} have constructed bijections between the primitive ideal
space of the quantized coordinate ring $\O_q(G)$ of a semisimple Lie group $G$
and the set of symplectic leaves in $G$ corresponding to a Poisson structure
which arises from the quantization process. In this ``standard'' case, the
symplectic leaves in $G$ are locally closed subvarieties, and they correspond
to the Poisson primitive ideals in the classical coordinate ring $\O(G)$.
Moreover, it is conjectured that the primitive ideal space of $\O_q(G)$ and the
Poisson primitive ideal space of $\O(G)$, with their respective Zariski
topologies, are homeomorphic. Results and conjectures such as these draw a focus
on the ideal theory of Poisson algebras. Our main concern in this paper is the
problem of identifying the Poisson primitive ideals in a Poisson algebra. To that
end, we establish a Poisson version of the famous Dixmier-Moeglin equivalence,
which applies to Poisson algebras equipped with suitable torus actions, a
hypothesis satisfied in many examples of interest. We also look at the Zariski
topology on the Poisson primitive spectrum of a Poisson algebra $R$, and prove
that it coincides with a natural quotient topology from the maximal ideal space
of $R$ in the case that the Poisson Dixmier-Moeglin equivalence holds.

In order to provide further detail, a few definitions are in order. A {\it
Poisson algebra\/} is a commutative algebra $R$ over a field $k$ (usually assumed
to have characteristic zero) together with an antisymmetric bilinear map
$\{-,-\}: R\times R \rightarrow R$ satisfying two key properties: the Jacobi
identity, so that $R$ together with $\{-,-\}$ forms a Lie algebra, and the
Leibniz rule, meaning that the maps $\{a,-\}$ (for $a\in R$) are derivations on
$R$. The ideals $I$ of $R$ for which $\{R,I\} \subseteq I$ are called {\it Poisson
ideals\/}, and within each ideal $J$ of $R$ there is a largest Poisson ideal,
which we call (following \cite{\BrGr}) the {\it Poisson core\/} of $J$. The {\it
Poisson primitive ideals\/} of $R$, finally, are the Poisson cores of the maximal
ideals. 

Recall that a noetherian algebra $A$ is said to satisfy the {\it Dixmier-Moeglin
equivalence\/} provided the following types of prime ideals in $A$ coincide: the
primitive ideals; the {\it locally closed\/} prime ideals, meaning those which
constitute locally closed points in the prime spectrum of $A$; and the {\it
rational\/} prime ideals, meaning those prime ideals $P$ in
$A$ such that the center of the Goldie quotient ring of $A/P$ is algebraic over
the base field. This equivalence was first established by Dixmier \cite{\Dix}
and Moeglin \cite{\Moe} for the enveloping algebra $U(\gfrak)$ of any finite
dimensional complex Lie algebra $\gfrak$. Among other settings where the
equivalence has been established is a class of algebras equipped with torus
actions studied by Letzter and the author \cite{\specstrat}. Our present
work uses many ideas from \cite{\specstrat}, but we follow the route laid out in
\cite{\BrGd, Chapters II.1--II.3, II.7--II.8} (see \cite{\BrGd, Theorem II.8.4}
for the equivalence result).

Taking the natural Poisson analogs of the above ideas, one says that a Poisson
algebra $R$ satisfies the {\it Poisson Dixmier-Moeglin equivalence\/} provided the
Poisson primitive ideals of $R$ are precisely the locally closed points of the
prime Poisson spectrum, and also precisely those prime Poisson ideals $P$ such
that the Poisson center of the quotient field of $R/P$ is algebraic over the base
field. Some examples in which this equivalence holds have been given by Oh
\cite{\Oh, Theorem 2.4, Proposition 2.13}, and Brown and Gordon have shown that it
holds for affine Poisson algebras with only finitely many Poisson primitive ideals
\cite{\BrGr, Lemma 3.4}. Our primary goal here is to establish the Poisson
Dixmier-Moeglin equivalence for Poisson algebras $R$ for which there is an
algebraic torus
$H$, acting rationally on $R$ by Poisson automorphisms, such that $R$ has only
finitely many prime Poisson ideals stable under $H$. As in \cite{\specstrat}, we
obtain an additional equivalence in our main theorem, based on a stratification
of the prime Poisson spectrum of $R$ arising from the action of $H$. Namely, the
Poisson primitive ideals of $R$ are (under the given hypotheses) exactly those
prime Poisson ideals which are maximal in their $H$-strata. (See below for
precise definitions.) 

In any Poisson algebra $R$, the process of taking Poisson
cores defines a canonical surjection from the maximal ideal space, $\max R$, onto
the Poisson primitive spectrum, $\Pprim R$. We show that if $R$ satisfies the 
Poisson Dixmier-Moeglin equivalence, then the Zariski topology on $\Pprim R$
coincides with the quotient topology induced by the above surjection. In
particular, if $R$ is the coordinate ring of an affine variety $V$ over an
algebraically closed field,
$\Pprim R$ becomes a topological quotient of $V$. Combining this result
with our main theorem, we thus obtain a large class of Poisson algebras in which
the Poisson primitive spectrum is a topological quotient of the maximal ideal
space. The methods are taken from joint work with Letzter on the quotient topology
problem for primitive spectra of quantized algebras \cite{\qaffquo}.
\medskip

We give a minimal sketch of Poisson structures and symplectic leaves in Section 4.
For readers interested in further background on these and
related subjects, we mention the following small sample of the available
literature:
\cite{\BrGd, Chapter III.5}, \cite{\CaWe, Chapter 3}, \cite{\EtSc, Lecture 1},
\cite{\KoSo, Chapter 1}. For some algebraic approaches, see \cite{\Far},
\cite{\FaLe}.
\medskip

The Lie algebra structure provided by a Poisson bracket plays no role in our
proofs; in fact, all that is needed is the set of derivations $\{a,-\}$.
Consequently, our results can be stated and proved in the context of a
commutative algebra equipped with a set of derivations. We derive everything at
this level of generality, in the first three sections of the paper. In the final
section, we specialize to the case of Poisson algebras, and discuss various
examples, in some of which the Poisson primitive ideals correspond to the
symplectic leaves in a complex affine Poisson variety.

Let us fix a base field:
\roster
\item"" {\it Throughout the paper, $k$ will denote a field of characteristic
zero.}
\endroster

\head 1. Ideals stable under derivations \endhead

As noted in the Introduction, our basic object of study will be a commutative
algebra equipped with a set of derivations. Readers who wish to concentrate on
Poisson algebras could do so by first comparing the beginnings of this section and
Section 4, and then substituting the prefix ``Poisson'' for ``$\Delta$-'' in what
follows.

The following terminology and
notation will be convenient. Many of these concepts are standard, and the basic
properties recorded in Lemma 1.1 mostly hold without any assumption of
commutativity, but for consistency we shall impose commutativity throughout. 

\definition{Definitions} A {\it commutative differential $k$-algebra\/} is a pair
$(R,\Delta)$ where $R$ is a commutative $k$-algebra and $\Delta$ a set of
$k$-linear derivations on $R$. We make no assumptions about any structure on
$\Delta$ -- in particular, it need not be a Lie subalgebra nor even a linear
subspace of $\Der_k(R)$. The {\it $\Delta$-center\/} of $R$ is the set
$$Z_\Delta(R)= \{r\in R\mid \delta(r)=0 \text{\;for all\;} \delta\in\Delta
\},$$  
a $k$-subalgebra of $R$. 

A {\it $\Delta$-ideal\/} of $R$ is any ideal $I$ of $R$ that is stable under
$\Delta$, that is, $\delta(I)\subseteq I$ for all $\delta\in\Delta$. The
{\it $\Delta$-core\/} of an arbitrary ideal $J$ in $R$ is the largest
$\Delta$-ideal contained in $J$, which we denote $(J:\Delta)$. It may be
described as follows \cite{\Dixbook, Lemma 3.3.2}:
$$(J:\Delta)= \{r\in J\mid \delta_1\delta_2 \cdots\delta_n(r)\in J \text{\;for
all\;} \delta_1,\dots,\delta_n \in\Delta\}.$$
Recall that the concept of a $\Delta$-prime ideal is
obtained by replacing arbitrary ideals by $\Delta$-ideals in the definition of a
prime ideal: A {\it $\Delta$-prime ideal\/} of $R$ is any proper $\Delta$-ideal
$Q$ of $R$ such that whenever $I$ and $J$ are $\Delta$-ideals of $R$ with
$IJ\subseteq Q$, either $I\subseteq Q$ or $J\subseteq Q$. In particular,
$(P:\Delta)$ is $\Delta$-prime for all prime ideals $P$. 

Adapting terminology from the theory of Poisson algebras, we define a {\it
$\Delta$-primitive ideal\/} of
$R$ to be any $\Delta$-ideal of the form $(M:\Delta)$ where $M$ is a maximal
ideal of $R$. 

The {\it $\Delta$-prime spectrum\/} of $R$ is the set $\delspec R$
of all $\Delta$-prime ideals, equipped with the natural Zariski topology.
Similarly, the {\it $\Delta$-primitive spectrum\/} is the subset $\delprim R$
consisting of all $\Delta$-primitive ideals, also equipped with the Zariski
topology. 
\enddefinition

\proclaim{Lemma 1.1} Let $(R,\Delta)$ be a commutative differential $k$-algebra.

{\rm (a)} $(P:\Delta)$ is prime for all prime ideals
$P$ of $R$.

{\rm (b)} Every $\Delta$-primitive ideal of $R$ is prime.

{\rm (c)} Every prime ideal minimal over a $\Delta$-ideal of $R$ is a
$\Delta$-ideal.

{\rm (d)} If $R$ is noetherian, every $\Delta$-prime ideal of $R$ is prime.

{\rm (e)} If $R$ is affine over $k$, every $\Delta$-prime ideal of
$R$ is an intersection of $\Delta$-primitive ideals. \endproclaim

\demo{Proof} (a) \cite{\Dixbook, Lemma 3.3.2}.

(b)(c) These are immediate from (a).

(d) Let $Q$ be a $\Delta$-prime of $R$. There exist prime ideals $Q_1,\dots,Q_n$
minimal over $Q$ such that $Q_1Q_2\cdots Q_n\subseteq Q$. The $Q_i$ are
$\Delta$-ideals by (c), so the $\Delta$-primeness of $Q$ implies that some
$Q_i=Q$.

(e) Let $Q$ be a $\Delta$-prime of $R$; then $Q$ is prime by (d). By the
Nullstellensatz,
$Q$ is an intersection of maximal ideals $M_i$, and so $Q= \bigcap_i
(M_i:\Delta)$.
\qed\enddemo

In view of Lemma 1.1(a) and the remarks above,
$$\delprim R \subseteq \delspec R \subseteq \spec R$$
for any commutative noetherian differential $k$-algebra $(R,\Delta)$.
Just as primitive ideals are more difficult to identify purely
ideal-theoretically than prime ideals, $\Delta$-primitive ideals are less
accessible than $\Delta$-prime ideals. A natural approach is to try to find
the $\Delta$-prime ideals first (even though there are more of them), and then to
develop criteria to tell which $\Delta$-prime ideals are
$\Delta$-primitive. Topological criteria (involving the space $\delspec R$) and
algebraic criteria (involving quotients of $R$ modulo $\Delta$-prime ideals) are
both useful. The key properties are given in terms of the following concepts.

\definition{Definitions} A {\it locally closed point\/} in a topological space
$X$ is any point which is (relatively) closed in some neighborhood. Note that a
point $x\in X$ is locally closed if and only if $\{x\}$ is the intersection of an
open and a closed set; hence, $x$ is locally closed if and only if $\{x\}$ is
(relatively) open in its closure.

A $\Delta$-prime ideal $P$ in a commutative noetherian differential $k$-algebra
$(R,\Delta)$ is {\it $\Delta$-ra\-tion\-al\/} provided the field $Z_\Delta(\Fract
R/P)$ is algebraic over $k$, where $\Fract R/P$ denotes the quotient field of
$R/P$ (recall from Lemma 1.1(d) that $R/P$ is a domain).
\enddefinition

The basic relations between these concepts and $\Delta$-primitivity were given in
the Poisson case by Oh \cite{\Oh, Propositions 1.7, 1.10}. Since the
proofs are slightly simpler in our case, and there is a change of
notation in the generalization to $(R,\Delta)$, we provide a sketch for the
reader's convenience.

\proclaim{Proposition 1.2} Let
$(R,\Delta)$ be a commutative differential
$k$-algebra, and assume that $R$ is affine over $k$. Then every locally closed
point in $\delspec R$ is $\Delta$-primitive, and every
$\Delta$-primitive ideal of $R$ is $\Delta$-rational.
\endproclaim

\demo{Proof} If $P$ is a locally closed point in $\delspec R$, there exist ideals
$I$ and $J$ in $R$ such that
$$\{P\}= \{Q\in \delspec R\mid Q\supseteq I\} \cap \{Q\in\delspec R\mid Q
\nsupseteq J\},$$
from which we see that
$$\bigcap\, \{Q\in\delspec R\mid Q \supsetneq P\} \supseteq P+J \supsetneq P.$$
By Lemma 1.1(e), $P$ is an intersection of
$\Delta$-primitive ideals, each of which is $\Delta$-prime. One of these must
equal $P$, proving that $P$ is $\Delta$-primitive.

Now let $P$ be an arbitrary $\Delta$-primitive ideal of $R$, and write
$P=(M:\Delta)$ for some maximal ideal $M$ of $R$. Since $R/M$ is finite
dimensional over $k$, it suffices to show that the $\Delta$-center of
the field $F= \Fract R/P$ embeds in $R/M$. After replacing $R$ by $R/P$, there
is no loss of generality in assuming that $P=0$.

We claim that $Z_\Delta(F)$ is contained in the localization $R_M$. Given a
fraction $ab^{-1}\in Z_\Delta(F)$, we have
$$0= \delta(ab^{-1})= \bigl( \delta(a)b-a\delta(b) \bigr)b^{-2}$$
and so $a\delta(b)=\delta(a)b$, for any $\delta\in\Delta$. Hence, $ab^{-1}=
\delta(a)\delta(b)^{-1}$. Repeating this argument, we see that
$$ab^{-1}= \delta_1(a)\delta_1(b)^{-1}= \delta_2\delta_1(a)
\delta_2\delta_1(b)^{-1}= \cdots= \bigl[ \delta_n\cdots \delta_2\delta_1(a) \bigr]
\bigl[ \delta_n\cdots \delta_2\delta_1(b) \bigr]^{-1}$$
for any $\delta_1,\dots,\delta_n \in\Delta$. Since $b\ne0$ while
$(M:\Delta)=P=0$, there exist $\delta_1,\dots,\delta_n \in\Delta$ such that
$\delta_n\cdots \delta_2\delta_1(b) \notin M$, and the claim is proved.

Therefore we obtain a $k$-algebra homomorphism
$$\phi : Z_\Delta(F) @>{\;\subseteq\;}>> R_M @>{\;\text{quo}\;}>> R_M/MR_M
@>{\;\cong\;}>> R/M.$$ 
Since $Z_\Delta(F)$ is a field, $\phi$ is an embedding, as
desired.
\qed\enddemo

The topological relationship between the spaces of prime and $\Delta$-prime
ideals in a commutative noetherian differential algebra is given by the following
result, which is a special case of \cite{\qaffquo, Proposition 1.7(c)}. We
provide a proof in order to avoid setting up the machinery of \cite{\qaffquo}.
Recall that a {\it topological quotient\/} of a topological space $X$ is a space
$Y$ together with a surjection $\pi: X\rightarrow Y$ such that the topology on
$Y$ is the quotient topology induced by $\pi$, that is,  the closed subsets of
$Y$ are exactly those subsets $W\subseteq Y$ for which $\pi^{-1}(W)$ is closed in
$X$.

\proclaim{Theorem 1.3} Let $(R,\Delta)$ be a commutative noetherian differential
$k$-algebra. Then the rule $\pi(Q)= (Q:\Delta)$ defines a continuous retraction
$$\pi : \spec R\rightarrow \delspec R,$$
and $\delspec R$ is a topological quotient of $\spec R$ via $\pi$. \endproclaim

\demo{Proof} By Lemma 1.1(d), $\delspec R\subseteq \spec R$. Thus, the given
rule defines a set-theoretic retraction of $\spec R$ onto $\delspec R$.

Let $V$ be a closed subset of $\delspec R$, and write $V= \{P\in\delspec R\mid
P\supseteq I\}$ for some ideal $I$ of $R$. Since we may replace $I$ by the
intersection of the ideals in $V$, there is no loss of generality in assuming that
$I$ is a
$\Delta$-ideal. Hence, 
$$\pi^{-1}(V)= \{Q\in\spec R\mid (Q:\Delta)\supseteq I\}= \{Q\in\spec R \mid
Q\supseteq I\},$$
a closed set in $\spec R$. This shows that $\pi$ is continuous.

Now consider a set $W\subseteq \delspec R$ such that $\pi^{-1}(W)$ is closed
in $\spec R$, say
$$\pi^{-1}(W)= \{Q\in\spec R \mid
Q\supseteq I\}$$
for some ideal $I$ of $R$. Since $\pi(P)=P$ for $P\in\delspec R$, we have
$$W= \pi^{-1}(W) \cap \delspec R= \{P\in\delspec R\mid P\supseteq I\},$$
a closed set in $\delspec R$. Therefore the topology on $\delspec R$ is indeed
the quotient topology induced by $\pi$.
\qed\enddemo

In the situation of Theorem 1.3, the map $\pi$ sends $\max R$ to
$\delprim R$. This restriction is continuous because $\pi$ is, and it is
surjective by definition of $\delprim R$. However, it need not be a
topological quotient map. For example, let $R= k[x]_{(x^2-x)}$, the localization
of a polynomial ring $k[x]$ at the semimaximal ideal $(x^2-x)$, and let
$\Delta= \{x\tfrac{d}{dx}\}$. Then $\delprim R$ consists of the two ideals
$(Rx:\Delta)=Rx$ and $(R(x-1):\Delta)=0$. The map $M\mapsto (M:\Delta)$ from
$\max R$ to $\delprim R$ is a continuous bijection, but $\delprim R$ is not a
topological quotient of $\max R$, since $\max R$ has the discrete topology while
$\delprim R$ does not. A more satisfactory example would be one affine over $k$,
but this is an open problem:

\proclaim{Question 1.4} Let $(R,\Delta)$ be a commutative differential
$k$-algebra, with $R$ affine over $k$. Is the Zariski topology on $\delprim R$
equal to the quotient topology induced from the continuous surjection $\max
R\rightarrow \delprim R$ given by $M\mapsto (M:\Delta)$? \endproclaim

In order for $\delprim R$ to be a topological quotient of $\max R$, it suffices,
by \cite{\qaffquo, Proposition 1.8}, to have 
$$P= \bigcap\, \bigl[ \pi^{-1}(\{P\})\cap \max R \bigr]= \bigcap\, \{M\in\max
R\mid (M:\Delta)=P\}  \tag1.1$$
for all $P\in\delprim R$. A sufficient condition for (1.1), in turn, would be a
converse of Proposition 1.2, as shown in the proof below.

\proclaim{Theorem 1.5} Let
$(R,\Delta)$ be a commutative differential
$k$-algebra, and assume that $R$ is affine over $k$. Then the rule
$\pi_m(Q)= (Q:\Delta)$ defines a continuous surjection
$$\pi_m: \max R \rightarrow \delprim R.$$
If, in addition, every $\Delta$-primitive ideal of $R$ is locally closed in
$\delspec R$, then $\delprim R$ is a topological quotient of $\max R$ via $\pi_m$.
\endproclaim

\demo{Proof} We have already noted that $\pi_m$ is a continuous surjection. Now
assume that every $\Delta$-primitive ideal of $R$ is locally closed in
$\delspec R$.

We claim that (1.1) holds for any $P\in\delprim R$. By assumption, $P$ is
locally closed in $\delspec R$, and so we see (as in the proof of Proposition
1.2) that
$$\bigcap\, \{Q\in\delspec R\mid Q \supsetneq P\}=: L \supsetneq P.$$
Separate the maximal ideals containing $P$ as follows:
$$\M_1= \{M\in\max R\mid M\supseteq L\} \qquad\text{and}\qquad \M_2= \{M\in
\max R\mid M\supseteq P \text{\;but\;} M\not\supseteq L\};$$
then $P= \left(\bigcap\M_1\right) \cap \left(\bigcap\M_2\right)$ by the
Nullstellensatz. Since $P$ is prime and $\bigcap\M_1 \supseteq
L \supsetneq P$, we must have $\bigcap\M_2=P$. For $M\in\M_2$, the ideal
$(M:\Delta)$ is a $\Delta$-prime which contains $P$ but not $L$, whence
$(M:\Delta)=P$. This establishes (1.1). 

Now consider a set $W\subseteq \delprim R$ such that $\pi_m^{-1}(W)$ is closed
in $\max R$, say
$$\pi_m^{-1}(W)= \{M\in\max R \mid M\supseteq I\}$$
for some ideal $I$. If $P\in\delprim R$ and $P\supseteq I$, then $P=\pi_m(M)$
for some $M\in\max R$, whence $M\supseteq (M:\Delta)=P\supseteq I$ and so
$M\in \pi_m^{-1}(W)$, yielding $P=\pi_m(M)\in W$. Conversely, if $P\in W$, then
$\pi_m^{-1}(\{P\})\subseteq \pi_m^{-1}(W)$ and so every ideal in
$\pi_m^{-1}(\{P\})$ contains $I$. By (1.1), the intersection of the ideals in
$\pi_m^{-1}(\{P\})$ equals
$P$, whence $P\supseteq I$. Thus
$$W= \{P\in\delprim R\mid P\supseteq I\},$$
a closed set in $\delprim R$. Therefore the topology on $\delprim R$ is indeed
the quotient topology induced by $\pi_m$.
\qed\enddemo

\head 2. Graded fields \endhead

Our main results will involve commutative differential algebras graded by (free
abelian) groups, and certain localizations of these algebras will be what are
known as ``graded fields'' -- this refers to the graded analog of the concept of
a field (see below), rather than to the idea of a field equipped with a grading.
Here we develop some properties of commutative differential graded fields. In
particular, we show that in such an algebra the $\Delta$-prime ideals correspond
precisely to the prime ideals of the $\Delta$-center. 

\definition{Definitions} Let $G$ be a group, which we assume to be
abelian and written multiplicatively. Recall that a {\it $G$-graded
$k$-algebra\/} is a $k$-algebra
$R$ equipped with a vector space direct sum decomposition $R= \bigoplus_{g\in G}
R_g$ such that $1\in R_1$ and $R_gR_h \subseteq R_{gh}$ for all $g,h\in G$. The
algebra $R$ is called {\it strongly $G$-graded\/} if $R_gR_h = R_{gh}$
for all $g,h\in G$. The subspaces $R_g$ are called the {\it homogeneous
components\/} of $R$ (with respect to the given grading), and the elements of a
given component $R_g$ are said to be {\it homogeneous of degree $g$\/}. Finally,
$R$ is called a {\it graded field\/} provided $R$ is a domain and all its nonzero
homogeneous elements are invertible. In that case, $R$ is necessarily strongly
graded.

A linear map $f:R\rightarrow R$ is {\it homogeneous of degree $d$\/} for some
$d\in G$ provided $f(R_g) \subseteq R_{gd}$ for all $g\in G$. 

To avoid repeating lengthy hypotheses, we define a {\it $G$-graded differential
$k$-algebra\/} to be a pair $(R,\Delta)$ such that
\roster
\item $R$ is a $G$-graded $k$-algebra;
\item $\Delta$ is a linear subspace of $\Der_k(R)$;
\item $\Delta= \bigoplus_{d\in G} \Delta_d$ where each $\Delta_d$ consists of
homogeneous derivations of degree $d$.
\endroster
(This concept is quite different than that of a ``differential graded'' algebra.)
\enddefinition

The simplest example of a $G$-graded field is the group algebra $kG$, equipped
with its standard grading, for any torsionfree abelian group $G$. (The nonzero
homogeneous elements of $kG$ are always invertible, but $kG$ is not a domain if
$G$ has torsion.) Versions of the results in this section have been proved for
Poisson structures on $kG$, where $G$ is free abelian of finite rank, by Oh and
Park \cite{\OhPa, Lemma 2.2, Theorem 2.3} and Vancliff \cite{\Van, Lemma 1.2}. 

\proclaim{Lemma 2.1} Let $G$ be an abelian group and $(R,\Delta)$ a
commutative $G$-graded differential
$k$-algebra. Assume that $R$ is a graded field.

{\rm (a)} The $\Delta$-center $Z_\Delta(R)$ is a
homogeneous subalgebra of
$R$, strongly graded by the subgroup $G_Z = \{x\in G \mid Z_\Delta(R)\cap R_x
\ne 0\}$ of $G$.

{\rm (b)} As $Z_\Delta(R)$-modules, $Z_\Delta(R)$ is a direct summand of $R$.

{\rm (c)} Suppose that $G_Z$ is free abelian of finite rank, with
basis $\{g_1,\dots ,g_n\}$. Choose a nonzero element $z_j\in Z_\Delta(R) \cap
R_{g_j}$ for each $j$. Then $Z_\Delta(R)$ is a Laurent polynomial ring of the form
$$Z_\Delta(R)= Z_\Delta(R_1) [z_1^{\pm1}, \dots, z_n^{\pm1}],$$
where the coefficient ring $Z_\Delta(R_1)= Z_\Delta(R) \cap R_1$ is a field. 
\endproclaim

\demo{Proof} (a) Since a homogeneous derivation $\delta$ on $R$ maps the
homogeneous components of an element $r\in R$ into distinct homogeneous components
of $R$, we see that $\delta$ cannot vanish on $r$ unless it vanishes on the
components of $r$. Hence, $Z_\Delta(R)$, which equals the intersection of the
kernels of the homogeneous derivations in $\Delta$, must be a homogeneous
subalgebra of $R$. Consequently, $Z_\Delta(R)= \bigoplus_{x\in G_Z}
Z_\Delta(R) \cap R_x$ where $G_Z$ is the subset of $G$ defined above. Since $R$
is a graded field, it is clear that $G_Z$ must be a subgroup of $G$.

(b) and (c) are proved exactly as \cite{\specstrat, Lemma 6.3} (cf\. \cite{\BrGd,
Lemma II.3.7}).
\qed\enddemo

\proclaim{Proposition 2.2} Let $G$ be an abelian group and $(R,\Delta)$ a
commutative $G$-graded differential
$k$-algebra. Assume that $R$ is a graded field. Then contraction
{\rm(}$I\mapsto I\cap Z_\Delta(R)$\/{\rm)} and extension {\rm(}$J\mapsto
RJ$\/{\rm)} provide inverse isomorphisms between the lattice of
$\Delta$-ideals of $R$ and the lattice of ideals of $Z_\Delta(R)$. \endproclaim

\demo{Proof} It suffices to prove that
\roster
\item"(a)" $I =
R(I\cap Z_\Delta(R))$ for any $\Delta$-ideal $I$ of $R$;
\item"(b)" $J =
(RJ)\cap Z_\Delta(R)$ for any ideal $J$ of $Z_\Delta(R)$.
\endroster
Statement (b) is clear from Lemma 2.1(b).

(a) Set $J= I\cap Z_\Delta(R)$, and suppose that $I\ne RJ$. Choose an element
$r\in I\setminus RJ$ with support $\{x_1,\dots, x_n\}$ of minimal
cardinality, and write $r= r_1+\dots +r_n$ for some nonzero elements
$r_i\in R_{x_i}$.

We claim that there cannot exist a nonzero element $s\in I$ whose support
is  contained in $\{x_2,\dots,x_n\}$. If such an element does exist, then after
possibly renumbering, we may assume that the support of $s$ includes $x_n$.
Write $s= s_2+\dots+ s_n$, with each $s_i\in R_{x_i}$ and $s_n\ne 0$. The
element $t = r_n^{-1}r -s_n^{-1}s$ is then an element in $I$
with support contained in $\{x_n^{-1}x_1,\dots, x_n^{-1}x_{n-1}\}$. By the
minimality of $n$, both $s$ and $t$ must lie in
$RJ$. But then $r\in RJ$, contradicting our
assumptions. This establishes the claim.

Now set $r'= r_1^{-1}r$, which is an element of $I$ with support
$\{1,x_1^{-1}x_2,\dots, x_1^{-1}x_n\}$ and identity component 1. Then $r'\in
I\setminus RJ$, and in particular $r'\notin Z_\Delta(R)$. Hence, there is some
$\delta\in \Delta$ such that $\delta(r')\ne 0$. Because of our hypotheses, we
may assume that $\delta$ is homogeneous of some degree $d\in G$. Since
$\delta(1)=0$, the support of $\delta(r')$ is contained in
$\{x_1^{-1}x_2d,\dots, x_1^{-1}x_nd\}$. There is some $i\ge2$ such that
$\delta(r_1^{-1}r_i)\ne0$, and then $u= r_i^{-1}\delta(r_1^{-1}r_i)$ is a
nonzero homogeneous element of degree $x_1^{-1}d$. But then $u^{-1}\delta(r')$
is a nonzero element of $I$ with support contained in $\{x_2,\dots,x_n\}$,
contradicting the claim above.

Therefore $I=RJ$. \qed\enddemo

\proclaim{Corollary 2.3} Let $G$ be an abelian group and $(R,\Delta)$ a
commutative $G$-graded differential
$k$-algebra. Assume that $R$ is a graded field. Then contraction and extension
provide mutually inverse homeomorphisms between $\delspec R$ and $\spec
Z_\Delta(R)$. Moreover, if $R$ is affine over $k$, contraction and extension
provide mutually inverse homeomorphisms between $\delprim R$ and $\max
Z_\Delta(R)$.
\endproclaim

\demo{Proof} Suppose that contraction maps $\delspec R$ to $\spec
Z_\Delta(R)$, and that extension maps $\spec Z_\Delta(R)$ to $\delspec R$. Then
Proposition 2.2 implies that these restricted maps are mutually inverse
bijections. Further, since contraction and extension preserve inclusions, the
proposition shows that the restricted maps are both continuous, hence
homeomorphisms.

Consider $Q\in\delspec R$. If $I$ and $J$ are ideals of
$Z_\Delta(R)$ such that $IJ\subseteq Q\cap Z_\Delta(R)$, then $RI$ and $RJ$
are $\Delta$-ideals of $R$ such that $(RI)(RJ)\subseteq Q$. Hence,
$RI\subseteq Q$ or $RJ\subseteq Q$, and so $I\subseteq Q\cap
Z_\Delta(R)$ or $J\subseteq Q\cap Z_\Delta(R)$. This shows that $Q\cap
Z_\Delta(R)$ lies in $\spec Z_\Delta(R)$.

Conversely, consider $P\in \spec Z_\Delta(R)$. If $I$ and $J$ are $\Delta$-ideals
of $R$ such that $IJ\subseteq RP$, then by Proposition 2.2, $I=RI'$ and $J=RJ'$
where $I'= I\cap Z_\Delta(R)$ and $J'= J\cap Z_\Delta(R)$. Then $I'J'\subseteq
RP\cap Z_\Delta(R)= P$, whence $I'\subseteq P$ or $J'\subseteq P$, and
so $I\subseteq RP$ or $J\subseteq RP$. This shows that $RP\in \delspec R$.

Therefore contraction and extension do map $\delspec R$ and $\spec
Z_\Delta(R)$ into each other, as desired.

Now assume that $R$ is affine over $k$.  As above, it will suffice to show that
contraction and extension map $\delprim R$ and $\max
Z_\Delta(R)$ into each other.

If $Q\in \delprim R$, then $Q=(M:\Delta)$ for
some maximal ideal $M$ of $R$. Since $R/M$ is finite
dimensional over $k$, the contraction $M\cap Z_\Delta(R)$ is a maximal ideal of
$Z_\Delta(R)$. By Proposition 2.2, $Q'= R(M\cap Z_\Delta(R))$ is a maximal proper
$\Delta$-ideal of $R$, and $Q'\cap Z_\Delta(R)= M\cap Z_\Delta(R)$. Since
$Q'\subseteq M$, we must have $Q'=(M:\Delta)=Q$. Thus $Q\cap Z_\Delta(R)=
M\cap Z_\Delta(R)$, a member of $\max Z_\Delta(R)$.

Finally, if $N\in\max Z_\Delta(R)$, then $RN$ is a maximal proper
$\Delta$-ideal of $R$. Hence, $RN=(M:\Delta)$ for any maximal ideal
$M\supseteq RN$, and so $RN\in\delprim R$.

Therefore contraction and extension do map $\delprim R$ and $\max
Z_\Delta(R)$ into each other, as desired.
\qed\enddemo

In the situations of interest to us, gradings on algebras arise from rational
actions of algebraic tori, as follows. For further detail, see \cite{\BrGd,
Chapter II.2}.

\definition{Definitions} An {\it algebraic torus over $k$\/} is a group of the
form $H= (\kx)^r$, where $r$ is a nonnegative integer, called the {\it rank\/} of
$H$. (Strictly speaking, $H$ should be called ``the group of $k$-rational points
of the algebraic group $({\kbar}{}^\times)^r$.'') The group $H$ is also an
algebraic group, based on its natural structure as an algebraic variety
(specifically, a locally closed subset of the affine space $k^r$). A {\it
character\/} of $H$ is any group homomorphism $H\rightarrow \kx$, and the {\it
rational characters\/} of $H$ are those which are also morphisms of algebraic
varieties. The set $\Hhat$ of rational characters of $H$ is a group under
pointwise multiplication of functions. (Note: Because $\Hhat$ is abelian, some
authors prefer to write it additively.) In fact,
$\Hhat$ is free abelian of rank $r$; a basis is given by the $r$ projection maps
$H= (\kx)^r \rightarrow \kx$. Having written $\Hhat$ as a multiplicative group,
questions of independence take an exponential form. Specifically, elements
$f_1,\dots,f_d\in \Hhat$ are {\it
$\ZZ$-linearly independent\/} if and only if the only integers $m_1,\dots,m_d$
for which $f_1^{m_1}f_2^{m_2} \cdots f_d^{m_d} = 1$ are $m_1= \cdots= m_d =0$.

Suppose that $H$ acts on a $k$-algebra $R$ by $k$-algebra automorphisms. A
nonzero element $r\in R$ is an {\it $H$-eigenvector\/} provided $h(r) \in kr$
for all $h\in H$. In that case, there is a unique character $f$ of $H$ such that
$h(r)= f(h)r$ for all $h\in H$, and $f$ is called the {\it $H$-eigenvalue\/} of
$r$. The collection of all $H$-eigenvectors in $R$ with $H$-eigenvalue $f$,
together with $0$, is a linear subspace called the {\it $H$-eigenspace of $R$
with $H$-eigenvalue $f$\/}. Note that the action of $H$ on $R$ induces an action
on $\Der_k(R)$ by $k$-linear automorphisms, where $h.\delta= h\circ\delta\circ
h^{-1}$ for $h\in H$ and $\delta\in \Der_k(R)$.

The action of $H$ on $R$ is {\it rational\/} provided 
\roster
\item The algebra $R$ is the direct sum of its $H$-eigenspaces (i.e., the action
of $H$ on $R$ is {\it semisimple\/});
\item The $H$-eigenvalues for the $H$-eigenspaces in $R$ are all rational.
\endroster
(This definition of a rational action is valid only for actions of algebraic
tori; see \cite{\BrGd, Definitions II.2.6} for the general
concept, and \cite{\BrGd, Theorem II.2.7} for the equivalence of the two
definitions.) When we have a rational action of $H$ on
$R$, we obtain a decomposition
$$R= \bigoplus_{f\in\Hhat} R_f$$
where $R_f$ denotes the $H$-eigenspace of $R$
with $H$-eigenvalue $f$. This decomposition turns $R$ into an $\Hhat$-graded
$k$-algebra (cf\. \cite{\BrGd, \S II.2.10}).

Finally, suppose that $(R,\Delta)$ is a differential $k$-algebra. We shall say
that {\it $H$ acts rationally on $(R,\Delta)$\/} provided
\roster
\item We are given a rational action of $H$ on $R$ by $k$-algebra automorphisms;
\item $\Delta$ is a linear subspace of $\Der_k(R)$, stable under the induced
$H$-action;
\item $\Delta$ is the direct sum of its $H$-eigenspaces, and the corresponding
$H$-eigenvalues are all rational.
\endroster
Observe that if $\delta\in\Delta$ is an $H$-eigenvector with $H$-eigenvalue $d$,
then $\delta$ is homogeneous of degree $d$ with respect to the $\Hhat$-grading on
$R$. Namely, since $h\delta h^{-1}= d(h)\delta$ for $h\in H$, we have
$$h\bigl( \delta(r) \bigr)= d(h)\delta \bigl( h(r) \bigr)= d(h)\delta \bigl(
f(h)r \bigr)= (df)(h) \delta(r)$$
for $f\in \Hhat$ and $r\in R_f$, whence $\delta(R_f) \subseteq R_{df}$.
Thus, when $H$ acts rationally on $(R,\Delta)$, the pair $(R,\Delta)$
becomes an $\Hhat$-graded differential $k$-algebra.
\enddefinition

\proclaim{Proposition 2.4} Let $(R,\Delta)$ be a commutative differential
$k$-algebra and
$H= (\kx)^r$ an algebraic torus acting rationally on $(R,\Delta)$. Assume that
$R$ is a graded field {\rm(}with respect to its $\Hhat$-grading\/{\rm)}.

{\rm (a)} The $\Delta$-center $Z_\Delta(R)$ is a Laurent polynomial ring, in at
most $r$ indeterminates, over the fixed field $Z_\Delta(R)^H$, which coincides
with the fixed field $Z_\Delta(\Fract R)^H$. The indeterminates can be chosen to
be $H$-eigenvectors with $\ZZ$-linearly independent $H$-eigenvalues.

{\rm (b)} Every $\Delta$-ideal of $R$ is generated by its intersection with
$Z_\Delta(R)$.

{\rm (c)} Contraction and extension
provide mutually inverse homeomorphisms between the spaces $\delspec R$ and $\spec
Z_\Delta(R)$. If $R$ is affine over $k$, then contraction and extension also
provide mutually inverse homeomorphisms between $\delprim R$ and $\max
Z_\Delta(R)$. \endproclaim

\demo{Proof} Except for the equality $Z_\Delta(R)^H= Z_\Delta(\Fract R)^H$, these
statements follow from Lemma 2.1, Proposition 2.2, and Corollary 2.3.

The inclusion $Z_\Delta(R)^H \subseteq Z_\Delta(\Fract R)^H$ is clear. Given an
element $u\in Z_\Delta(\Fract R)^H$, observe that the set $I= \{r\in R\mid ru\in
R\}$ is a nonzero $H$-stable ideal of $R$. Since $H$ acts semisimply on $R$, it
follows that $I$ is spanned by $H$-eigenvectors, that is, $I$ is a homogeneous
ideal (with respect to the $\Hhat$-grading on $R$). Since $R$ is a graded field,
$I=R$, whence $u\in R$. Therefore $u\in Z_\Delta(R)^H$, as desired. \qed\enddemo

\head 3. Stratification\endhead 

We now investigate the general case of a commutative differential
$k$-algebra equipped with a rational torus action. The results of the previous
section will be applied to suitable localizations of factor algebras. 

\definition{Definitions} Let $H$ be a group acting on a ring $R$ by
automorphisms. An {\it $H$-ideal\/} of $R$ is any ideal $I$ of $R$ that is stable
under $H$, that is, $h(I)=I$ for all $h\in H$. (It is sufficient to check that
$h(I)\subseteq I$ for all $h\in H$.) Given an arbitrary ideal $I$ in $R$, let
$(I:H)$ denote the largest $H$-ideal of $R$ contained in $I$, that is, $(I:H)=
\bigcap_{h\in H} h(I)$. An {\it $H$-prime ideal\/} of $R$ is any proper $H$-ideal
$Q$ of $R$ such that whenever $I$ and $J$ are $H$-ideals of $R$ with $IJ\subseteq
Q$, either
$I\subseteq Q$ or $J\subseteq Q$. In particular, $(P:H)$ is $H$-prime for any
prime ideal $P$.

Now suppose that we also have a set $\Delta$ of derivations on $R$. An
{\it $(H,\Delta)$-ideal\/} of $R$ is any ideal of $R$ that is stable under both
$H$ and
$\Delta$. We then define the notion of an {\it $(H,\Delta)$-prime ideal\/} in
parallel with $H$-prime or $\Delta$-prime ideals, and we write $\hdelspec R$ for
the {\it $(H,\Delta)$-spectrum\/} of $R$, that is, the set of all
$(H,\Delta)$-prime ideals of $R$, equipped with the natural Zariski topology.
\enddefinition

\proclaim{Lemma 3.1} Let $(R,\Delta)$ be a commutative noetherian differential
$k$-algebra and $H$ an algebraic torus acting rationally on
$(R,\Delta)$. 

{\rm (a)} If $I$ is a $\Delta$-ideal of $R$, then $(I:H)$ is an
$(H,\Delta)$-ideal. Similarly, if $I$ is an $H$-ideal of $R$, then $(I:\Delta)$
is an $(H,\Delta)$-ideal.

{\rm (b)} $(P:H) \in \hdelspec R$ for all $P\in \delspec R$.

{\rm (c)} The following sets coincide:

\qquad {\rm (1)} $\hdelspec R$;

\qquad {\rm (2)} The set of all $H$-prime $\Delta$-ideals in $R$;

\qquad {\rm (3)} The set of all $\Delta$-prime $H$-ideals in $R$;

\qquad {\rm (4)} The set of all prime $(H,\Delta)$-ideals in $R$.
\endproclaim

\demo{Proof} (a) Assume first that $I$ is a $\Delta$-ideal of $R$. If $\delta\in
\Delta$ is an $H$-eigenvector with $H$-eigenvalue $d$, then
$$h\delta(I:H)= d(h)\cdot\delta h(I:H)\subseteq \delta(I)\subseteq I$$
for all $h\in H$, whence $\delta(I:H)\subseteq (I:H)$. It follows that $(I:H)$ is
stable under $\Delta$, and so it is an $(H,\Delta)$-ideal.

Now assume that $I$ is an $H$-ideal of $R$, and consider $h\in H$. For any
$H$-eigenvectors $\delta_1,\dots,\delta_n\in \Delta$ with respective
$H$-eigenvalues $d_1,\dots,d_n$, we have
$$\delta_1\delta_2\cdots\delta_nh \bigl( (I:\Delta) \bigr)= \bigl[ (d_1d_2\cdots
d_n)(h) \bigr]^{-1} h\delta_1\delta_2\cdots\delta_n \bigl( (I:\Delta) \bigr)
\subseteq h(I)= I,$$ 
from which it follows that $h\bigl( (I:\Delta) \bigr) \subseteq (I:\Delta)$.
Hence, $(I:\Delta)$ is stable under $H$, and so it is an $(H,\Delta)$-ideal.

(b) If $P\in \delspec R$, then $P$ is prime by Lemma 1.1(d), and so $(P:H)$ is an
$H$-prime ideal. By part (a), $(P:H)$ is an $(H,\Delta)$-ideal, and we conclude
from the $H$-primeness of $(P:H)$ that it must be $(H,\Delta)$-prime.

(c) Since all $\Delta$-prime ideals of $R$ are prime (Lemma 1.1), as are all
$H$-prime ideals (e.g., \cite{\BrGd, Proposition II.2.9}), the sets (2), (3), and
(4) coincide. Clearly $(4) \subseteq (1)$.

Now let $Q$ be an $(H,\Delta)$-prime ideal of $R$. Any prime ideal minimal over
$Q$ is a $\Delta$-ideal by Lemma 1.1. Since $R$ is noetherian, there are only
finitely many prime ideals minimal over $Q$, and they are permuted by $H$, so
their $H$-orbits are finite. Thus, \cite{\BrGd, Proposition II.2.9} implies that
the prime ideals minimal over $Q$ are $H$-ideals, and so they all lie in the
set (4). Now there exist prime ideals $Q_1,\dots,Q_n$ minimal over $Q$ such that
$Q_1Q_2 \cdots Q_n \subseteq Q$. Since $Q$ is $(H,\Delta)$-prime, we conclude
that some $Q_j=Q$. Therefore $(1)\subseteq (4)$. \qed\enddemo

\definition{Definitions} Let $(R,\Delta)$ be a commutative noetherian differential
$k$-algebra and $H$ an algebraic torus acting rationally on
$(R,\Delta)$. For $J\in \hdelspec R$, we define the {\it $J$-stratum\/} of
$\delspec R$ to be the set
$$\delspec_J R = \{ P\in\delspec R \mid (P:H)=J \}.$$
In view of Lemma 3.1(b), we obtain a partition
$$\delspec R= \bigsqcup_{J\in\hdelspec R} \delspec_J R,$$
which we refer to as the {\it $H$-stratification of $\delspec R$\/}. We define
strata denoted $\delprim_J R$ in $\delprim R$ in the same way as $\delspec_J R$,
and obtain a corresponding
$H$-stratification of that space:
$$\delprim R= \bigsqcup_{J\in\hdelspec R} \delprim_J R.$$

For any $J\in\hdelspec R$, let $\E_J$ denote the set of
$H$-eigenvectors in $R/J$. Since $R/J$ is a domain (by Lemma 3.1), $\E_J$ is
multiplicatively closed, and the localization $R_J= (R/J)[\E_J^{-1}]$ is a
subalgebra of $\Fract (R/J)$. Note that the actions of $H$ and $\Delta$ on $R$
both extend naturally to $R/J$ and $R_J$, and then to $\Fract R/J$, so that we have
commutative differential $k$-algebras $(R/J,\Delta)$ and $(R_J,\Delta)$, as well
as $(\Fract R/J, \Delta)$.
\enddefinition

\proclaim{Theorem 3.2} Let $(R,\Delta)$ be a commutative noetherian differential
$k$-algebra, and $H= (\kx)^r$ an algebraic torus acting rationally on
$(R,\Delta)$. Let $J$ be a prime $(H,\Delta)$-ideal of $R$.

{\rm (a)} The algebra $R_J$ is a graded field, with respect to its induced
$\Hhat$-grading.

{\rm (b)} $\delspec_J R$ is homeomorphic to $\delspec R_J$ via localization and
contraction.

{\rm (c)} $\delspec R_J$ is homeomorphic to $\spec Z_\Delta(R_J)$ via contraction
and extension.

{\rm (d)} $Z_\Delta(R_J)$ is a Laurent polynomial ring, in at most $r$
indeterminates, over the fixed field $Z_\Delta(R_J)^H= Z_\Delta(\Fract R/J)^H$.
The indeterminates can be chosen to
be $H$-eigenvectors with $\ZZ$-linearly independent $H$-eigenvalues.
\endproclaim

\demo{Proof} (a) With respect to the $\Hhat$-grading, $R_J$ is obtained from the
domain $R/J$ by inverting all nonzero homogeneous elements. Consequently, $R_J$
is a graded field.

(b) Standard localization theory yields that localization and contraction give
mutually inverse homeomorphisms between the sets
$$X_J := \{ P\in\spec R \mid P\supseteq J \text{\;and\;} (P/J)\cap \E_J =
\varnothing\}$$ 
and $\spec R_J$. It is clear that these maps restrict to mutually inverse
homeomorphisms between $X_J\cap \delspec R$ and $\delspec R_J$, so it just
remains to show that $X_J\cap \delspec R= \delspec_J R$. 

If $P\in \delspec_J R$, then $(P/J:H)=0$ and so $P/J$ contains no nonzero
$H$-ideals of $R/J$. Hence, $P/J$ contains no $H$-eigenvectors, that is, $P\in
X_J$. 

Conversely, if $P\in X_J\cap \delspec R$, then $P\supseteq J$ but $P/J$ contains
no $H$-eigenvectors of $R/J$. Since all $H$-ideals of $R/J$ are spanned by their
$H$-eigenvectors, $P/J$ contains no nonzero $H$-ideals, that is, $(P/J:H)=0$.
Thus, $(P:H)=J$ and $P\in \delspec_J R$, as desired.

(c)(d) These follow from (a) together with Proposition 2.4. \qed\enddemo

We can now give the following Dixmier-Moeglin equivalence for commutative
differential algebras equipped with rational torus actions.

\proclaim{Theorem 3.3} Let $(R,\Delta)$ be a commutative differential
$k$-algebra, and $H= (\kx)^r$ an algebraic torus acting rationally on
$(R,\Delta)$. 

Assume that $R$ is affine over $k$, and that it has
only finitely many prime $(H,\Delta)$-ideals. Let $J$ be one of them. For
$P\in \delspec_J R$, the following conditions are equivalent:

{\rm (a)} $P$ is locally closed in $\delspec R$.

{\rm (b)} $P$ is $\Delta$-primitive.

{\rm (c)} $Z_\Delta(\Fract R/P)$ is algebraic over $k$.

{\rm (d)} $P$ is maximal in $\delspec_J R$. \endproclaim

\demo{Proof} (a)$\Longrightarrow$(b)$\Longrightarrow$(c) by Proposition 1.2.

(d)$\Longrightarrow$(a): If $J$ is maximal in $\hdelspec R$,
then every $\Delta$-prime containing $J$ is in $\delspec_J R$. In this case,
$P$ is maximal in $\delspec R$, and thus is trivially locally closed in
$\delspec R$.

Now suppose that $J$ is not maximal in $\hdelspec R$, and let
$J_1,\dots,J_n$ be the prime $(H,\Delta)$-ideals of $R$ that properly contain
$J$. Then the ideal $\Jhat= J_1\cap \dots\cap J_n$ must
properly contain $J$. But $J=(P:H)$, so this implies that $\Jhat
\not\subseteq P$. Any $\Delta$-prime ideal $Q$ which properly contains $P$ cannot
lie in $\delspec_J R$. Since $(Q:H)\supseteq (P:H)=J$, it follows that
$(Q:H)=J_i$ for some $i$, and so $Q\supseteq J_i\supseteq \Jhat$. Thus, the
$\Delta$-prime ideals properly containing $P$ all contain $\Jhat$, which shows
that
$$\{P\}= \{Q\in \delspec R \mid Q\supseteq P\} \cap \{Q\in \delspec R \mid Q
\nsupseteq \Jhat\}.$$
 Therefore $P$ is locally closed in $\delspec R$.

(c)$\Longrightarrow$(d): For this part of the proof, we pass to $R/J$, and so
we may assume that $J=0$. Thus, $R$ is now a domain.

By Theorem 3.2(b), $P$ is disjoint from $\E_J$ and $P$ induces a $\Delta$-prime
ideal $PR_J$ in $R_J$. Note that $R_J/PR_J$ is a localization
of $R/P$, and so $\Fract R_J/PR_J= \Fract R/P$. Hence,
$$Z_\Delta(\Fract R_J/PR_J)= Z_\Delta(\Fract R/P),$$
which is algebraic over $k$ by hypothesis.

Set $Q= PR_J\cap Z_\Delta(R_J)$, which is a prime ideal of $Z_\Delta(R_J)$ by
Theorem 3.2(c). Further, the natural embedding
$$Z_\Delta(R_J)/Q \rightarrow R_J/PR_J \rightarrow \Fract R_J/PR_J$$
maps $Z_\Delta(R_J)/Q$ into $Z_\Delta(\Fract R_J/PR_J)$, and so
$Z_\Delta(R_J)/Q$ must be algebraic over $k$. It follows that
$Z_\Delta(R_J)/Q$ is a field, whence $Q$ is a  maximal ideal of
$Z_\Delta(R_J)$. By Theorem 3.2(b)(c), $PR_J$ is maximal in $\delspec R_J$, and
therefore $P$ is maximal in $\delspec_J R$.
\qed\enddemo

\proclaim{Corollary 3.4} Let $(R,\Delta)$ be a commutative differential
$k$-algebra, and $H= (\kx)^r$ an algebraic torus acting rationally on
$(R,\Delta)$. Assume that $R$ is affine over $k$, and that it has
only finitely many prime $(H,\Delta)$-ideals. Then the rule
$\pi_m(Q)= (Q:\Delta)$ defines a continuous surjection
$$\pi_m: \max R \rightarrow \delprim R,$$
and $\delprim R$ is a topological quotient of $\max R$ via $\pi_m$. \endproclaim

\demo{Proof} Theorems  3.3 and 1.5. \qed\enddemo

Finally, we tighten up the picture for the case that $k$ is algebraically closed.

\proclaim{Theorem 3.5} Let $(R,\Delta)$ be a commutative differential
$k$-algebra, and $H= (\kx)^r$ an algebraic torus acting rationally on
$(R,\Delta)$. Assume that $k$ is algebraically closed, that $R$ is affine over
$k$, and that $R$ has only finitely many prime $(H,\Delta)$-ideals.

{\rm (a)} For each prime $(H,\Delta)$-ideal $J$ of $R$, the algebra
$Z_\Delta(R_J)$ is a Laurent polynomial ring of the form
$k[z_1^{\pm1},\dots,z_n^{\pm1}]$, for some nonnegative integer $n=n(J) \le r$.
Consequently, the $\Delta$-primitive ideals within $\delspec_J R$ are precisely
the inverse images in $R$ of the ideals
$$(R/J) \cap \bigl( R_J(z_1-\alpha_1) +\cdots+ R_J(z_n-\alpha_n \bigr)
\vartriangleleft R/J,$$ 
for arbitrary nonzero scalars $\alpha_1,\dots,\alpha_n \in \kx$.

{\rm (b)} The $H$-orbits within $\delprim R$ coincide with the $H$-strata of
$\delprim R$. In particular, there are only finitely many $H$-orbits in
$\delprim R$.
\endproclaim

\demo{Proof} (a) In view of Theorems 3.2 and 3.3, the $\Delta$-primitive ideals
within $\delspec_J R$ are precisely the inverse images in $R$ of the ideals in
$R/J$ extended from maximal ideals of $Z_\Delta(R_J)$. Thus, we just need to show
that $Z_\Delta(R_J)$ is a Laurent polynomial algebra over $k$ in at most $r$
variables. That will follow from Theorem 3.2(d) once we show that
$Z_\Delta(R_J)^H=k$.

Since $\delspec_J R$ is nonempty (it contains $J$), there exists
a $\Delta$-prime ideal $P$ which is maximal in $\delspec_J R$. By Theorem 3.3,
the field $Z_\Delta(\Fract R/P)$ is algebraic over $k$, and thus equals $k$.
Now since $Z_\Delta(R_J)^H$ is a field, the natural homomorphism
$$Z_\Delta(R_J)^H @>{\;\subseteq\;}>>  Z_\Delta(R_J) \rightarrow
Z_\Delta(R_J/PR_J) @>{\;\subseteq\;}>>  Z_\Delta(\Fract R_J/PR_J) \cong
Z_\Delta(\Fract R/P)$$  
is an embedding. Therefore $Z_\Delta(R_J)^H=k$, as required.

(b) It is clear that the $H$-action on $\delprim R$ preserves the $H$-strata.
Hence, any $H$-orbit within $\delprim R$ is contained in a stratum $\delprim_J
R$, for some prime $(H,\Delta)$-ideal $J$ of $R$, and it only remains to show
that $H$ acts transitively on $\delprim_J R$.

Thus, let $P_1,P_2\in \delprim_J R$. By Theorems 3.2 and 3.3, the $P_i$ induce
maximal $\Delta$-ideals in $R_J$, which contract to maximal ideals $Q_i :=
P_iR_J\cap Z_\Delta(R_J)$ in $Z_\Delta(R_J)$. Now by Theorem 3.2(d),
we can choose the indeterminates $z_i$ in part (a) to be
$H$-eigenvectors with $\ZZ$-linearly independent $H$-eigenvalues $f_i$. Each
$$Q_i= \langle z_1-\alpha_{i1},\, \dots,\, z_n-\alpha_{in} \rangle$$
for some $\alpha_{ij}\in \kx$. Since $f_1,\dots,f_n$ are $\ZZ$-linearly
independent elements of $\Hhat$, there exists $h\in H$ such that $f_j(h)=
\alpha_{1j}\alpha_{2j}^{-1}$ for all $j$ (e.g., \cite{\Hum, Lemma 16.2C}).
Then
$$h(z_j- \alpha_{1j})= f_j(h)z_j- \alpha_{1j}= \alpha_{1j}\alpha_{2j}^{-1} (z_j-
\alpha_{2j})$$
for all $j$, whence $h(Q_1)= Q_2$. As a result,
$$h(P_1)R_J \cap Z_\Delta(R_J)= h(Q_1)= Q_2= P_2R_J \cap Z_\Delta(R_J),$$
and therefore we conclude from Theorem 3.2 that $h(P_1)= P_2$, as desired.
\qed\enddemo

\head 4. The Poisson case\endhead

\definition{Definitions} A {\it Poisson $k$-algebra\/} is a pair $(R,\{-,-\})$
where $R$ is a commutative $k$-algebra and $\{-,-\}$ is a {\it Poisson bracket\/}
on $R$, that is, a bilinear map $\{-,-\}: R\times R\rightarrow R$ such that
\roster
\item"(a)" The vector space $R$ equipped with the binary operation $\{-,-\}$ is a
Lie algebra over $k$, and
\item"(b)" For each $a\in R$, the $k$-linear map $\{a,-\} : R\rightarrow R$ is a
derivation (called the {\it Hamiltonian associated to $a$\/}).
\endroster
We typically assume that a Poisson bracket has been given and is denoted by
$\{-,-\}$, and hence will refer to $R$ itself as a Poisson algebra. A {\it
Poisson automorphism\/} of $R$ is any $k$-algebra automorphism
of $R$ which preserves the Poisson bracket.

All of the concepts defined at the beginning of Section 1 are considered for a
Poisson algebra $R$ relative to the set $\{R,-\}$ of Hamiltonians on $R$. Thus,
the {\it Poisson center\/} of $R$ is the $\{R,-\}$-center, which we shall denote
$Z_P(R)$. (The elements of the Poisson center are sometimes called {\it
Casimirs\/}, in which case $Z_P(R)$ is denoted $\operatorname{Cas} R$.) A {\it
Poisson ideal\/} of $R$ is any
$\{R,-\}$-ideal. Given an arbitrary ideal $J$ in $R$, there is a largest Poisson
ideal contained in $J$, which in the notation of Section 1 would be written
$(J:\{R,-\})$. Following \cite{\BrGr}, we call $(J:\{R,-\})$ the {\it Poisson
core of $J$\/}; we shall denote it $\Pcore(J)$. The {\it Poisson primitive
ideals\/} of
$R$ are the $\{R,-\}$-primitive ideals, that is, the Poisson cores of the maximal
ideals of $R$ (in
\cite{\Oh}, these are called {\it symplectic ideals\/}). A {\it Poisson-prime\/}
ideal of $R$ is any $\{R,-\}$-prime ideal. If $R$ is noetherian, then by Lemma
1.1(d), the Poisson-prime ideals in $R$ are precisely the ideals which are both
Poisson ideals and prime ideals; in that case, the hyphen in the term
``Poisson-prime'' becomes unnecessary. Finally, we write
$\Pspec R$ and
$\Pprim R$ for $\{R,-\}\operatorname{-spec} R$ and $\{R,-\}\operatorname{-prim}
R$, respectively.

In order to apply the results of previous sections to the Poisson setting, one
observation is needed: If $R$ is a Poisson algebra and $H$ is an algebraic torus
acting rationally on $R$ by Poisson automorphisms, then $H$ acts rationally on
$(R,\{R,-\})$. To see this, observe that given $h\in H$ and $a\in
R$, we have $h(\{a,h^{-1}(b)\})= \{h(a),b\}$ for all $b\in R$, whence $h.\{a,-\}=
h\circ\{a,-\}\circ h^{-1}= \{h(a),-\}$. In particular, if $a$ is an
$H$-eigenvector, then so is $\{a,-\}$, with the same eigenvalue. As every element
of $R$ is a sum of $H$-eigenvectors with rational $H$-eigenvalues, the same holds
in $\{R,-\}$. Therefore $H$ acts rationally on
$(R,\{R,-\})$, as claimed.
\enddefinition

\proclaim{Theorem 4.1} Let $R$ be a noetherian Poisson $k$-algebra.

{\rm (a)} The rule $\pi(Q)= \Pcore(Q)$ defines a continuous retraction
$$\pi : \spec R\rightarrow \Pspec R,$$
and $\Pspec R$ is a topological quotient of $\spec R$ via $\pi$.

{\rm (b)} Now suppose that $R$ is an affine $k$-algebra. Let $H=(k^\times)^r$ be
an algebraic torus acting rationally on $R$ by Poisson automorphisms, and
assume that
$R$ has only finitely many prime Poisson $H$-ideals. Then $\pi$ restricts to a
continuous surjection
$$\pi_m: \max R \rightarrow \Pprim R,$$
and $\Pprim R$ is a topological quotient of $\max R$ via $\pi_m$. \endproclaim

\demo{Proof} Theorem 1.3 and Corollary 3.4.
\qed\enddemo

As in Section 3, when we have a group $H$ acting on a Poisson algebra $R$ by
Poisson automorphisms, the Poisson spectrum $\Pspec R$ obtains an
{\it $H$-stratification\/}, namely the partition 
$$\Pspec R= \bigsqcup_{\text{prime Poisson $H$-ideals\ }J} \Pspec_J R,$$
where the {\it $H$-strata\/} $\Pspec_J R$ are given by
$$\Pspec_J R= \{ P\in \Pspec R \mid (P:H)= J\}.$$
There is a parallel $H$-stratification of $\Pprim R$.
Given a prime Poison $H$-ideal $J$ in $R$, we again denote the localization of
$R/J$ with respect to the multiplicative set of its $H$-eigenvectors by $R_J$,
and we equip
$R_J$ and
$R/J$, as well as
$\Fract R/J$, with the induced Poisson structures.  In Poisson notation, the
stratification theorem takes the following form:

\proclaim{Theorem 4.2} Let $R$ be a noetherian Poisson $k$-algebra, and $H=
(\kx)^r$ an algebraic torus acting rationally on
$R$ by Poisson automorphisms. Let $J$ be a prime Poisson $H$-ideal of $R$.

{\rm (a)} The algebra $R_J$ is a graded field, with respect to its induced
$\Hhat$-grading.

{\rm (b)} $\Pspec_J R$ is homeomorphic to $\Pspec R_J$ via localization and
contraction.

{\rm (c)} $\Pspec R_J$ is homeomorphic to $\spec Z_P(R_J)$ via contraction
and extension.

{\rm (d)} $Z_P(R_J)$ is a Laurent polynomial ring, in at most $r$
indeterminates, over the fixed field $Z_P(R_J)^H= Z_P(\Fract R/J)^H$.
The indeterminates can be chosen to
be $H$-eigenvectors with $\ZZ$-linearly independent $H$-eigenvalues.
\endproclaim

\demo{Proof} Theorem 3.2. \qed\enddemo

Our main theorem, given next, provides a Dixmier-Moeglin equivalence for Poisson
algebras equipped with rational torus actions.

\proclaim{Theorem 4.3} Let $R$ be an affine Poisson $k$-algebra, and $H=
(\kx)^r$ an algebraic torus acting rationally on $R$ by Poisson automorphisms.
Assume that $R$ has
only finitely many prime Poisson $H$-ideals, and let $J$ be one of them. For
$P\in \Pspec_J R$, the following conditions are equivalent:

{\rm (a)} $P$ is locally closed in $\Pspec R$.

{\rm (b)} $P$ is Poisson primitive.

{\rm (c)} $Z_P(\Fract R/P)$ is algebraic over $k$.

{\rm (d)} $P$ is maximal in $\Pspec_J R$. \endproclaim

\demo{Proof} Theorem 3.3. \qed\enddemo

Theorem 4.3 sets up a general framework which covers various previous examples
of the Poisson Dixmier-Moeglin equivalence, such as those of Oh \cite{\Oh,
Theorem 2.4, Proposition 2.13}. See below for further detail about these
examples. Note that if $R$ is an affine Poisson algebra with only finitely many
Poisson primitive ideals, then, by Lemma 1.1(e), $R$ has only finitely many prime
Poisson ideals. Thus, the case of Theorem 4.3 in which $H= \langle 1\rangle$
covers \cite{\BrGr, Lemma 3.4}.

\proclaim{Theorem 4.4} Let $R$ be an affine Poisson $k$-algebra, and $H=
(\kx)^r$ an algebraic torus acting rationally on $R$ by Poisson automorphisms.
Assume that $k$ is algebraically closed, and that $R$ has
only finitely many prime Poisson $H$-ideals.

{\rm (a)} For each prime Poisson $H$-ideal $J$ of $R$, the algebra 
$Z_P(R_J)$ is a Laurent polynomial ring of the form
$k[z_1^{\pm1},\dots,z_n^{\pm1}]$, for some nonnegative integer $n=n(J) \le r$.
Consequently, the Poisson primitive ideals within $\Pspec_J R$ are precisely
the inverse images in $R$ of the ideals
$$(R/J) \cap \bigl( R_J(z_1-\alpha_1) +\cdots+ R_J(z_n-\alpha_n \bigr)
\vartriangleleft R/J,$$ 
for arbitrary nonzero scalars $\alpha_1,\dots,\alpha_n \in \kx$.

{\rm (b)} The $H$-orbits within $\Pprim R$ coincide with the $H$-strata of
$\Pprim R$. In particular, there are only finitely many $H$-orbits in
$\Pprim R$.
\endproclaim

\demo{Proof} Theorem 3.5. \qed\enddemo

In the setting of Theorem 4.4, part (a) provides an explicit recipe for writing
down all the Poisson primitive ideals of $R$, provided one can find all the prime
Poisson $H$-ideals $J$ in $R$ and one can compute indeterminates for the Laurent
polynomial rings $Z_P(R_J)$. The following examples exhibit some
uses of this recipe.

\definition{Example 4.5} Let $R= k[x_1^{\pm1},\dots,x_n^{\pm1}]$ be a Laurent
polynomial algebra over $k$, and let $\pi= (\pi_{ij})$ be an $n\times n$
antisymmetric matrix over $k$. (At this point, $k$ does not need to be
algebraically closed.) There is a unique Poisson bracket on
$R$ such that
$$\{x_i,x_j\}= \pi_{ij}x_ix_j$$
for all $i,j$; a complete formula for this bracket is
$$\{f,g\}= \sum_{i,j=1}^n \pi_{ij}x_ix_j \frac{\partial f}{\partial x_i}
\frac{\partial g}{\partial x_j}  \tag4.1$$
for $f,g\in R$. The torus $H= (\kx)^n$ acts on $R$ by $k$-algebra automorphisms
such that
$$(\alpha_1,\dots,\alpha_n).x_i= \alpha_ix_i  \tag4.2$$
for $(\alpha_1,\dots,\alpha_n) \in H$ and $i=1,\dots,n$. Observe that this is a
rational action by Poisson automorphisms. It is easily checked that $0$ is the
only prime $H$-ideal of $R$, and thus the only prime Poisson
$H$-ideal. Hence, Theorem 4.3 implies that $R$ satisfies the Poisson
Dixmier-Moeglin equivalence, recovering a result of Oh \cite{\Oh, Theorem 2.4}.
Further, Theorem 4.1 shows that $\Pprim R$ is a topological quotient of $\max
R$, via the map $M\mapsto \Pcore(M)$. This result is implicit in the work of Oh,
Park, and Shin \cite{\OPS}. 

Since $\Pspec R$ consists of a single $H$-stratum, Theorem 4.3 also implies that
the Poisson primitive ideals of $R$ are precisely the maximal Poisson ideals. We
can get a handle on these with Theorem 4.2. Observe that the $H$-eigenvectors in
$R$ are just the monomials $x_1^{m_1}x_2^{m_2} \cdots x_n^{m_n}$, which are
already invertible in $R$. Thus, the algebra $R_0$ of Theorem 4.2 (corresponding
to the prime Poisson $H$-ideal $0$) is just $R$ itself. Consequently, the
theorem shows that the maximal Poisson ideals of $R$ are precisely the ideals
extended from maximal ideals of the Poisson center $Z_P(R)$. By \cite{\OhPa,
Lemma 2.2} or \cite{\Van, Lemma 1.2(a)}, 
$$Z_P(R) = k\text{-span of\ } \{ x_1^{m_1}x_2^{m_2} \cdots
x_n^{m_n} \mid m_i\in \ZZ \text{\ and\ } \sum_{i=1}^n m_i\pi_{ij} =0 \text{\ for
all\ } j \}.$$
Thus, $Z_P(R)$ can be a Laurent polynomial algebra in any number of
indeterminates from $0$ to $n$, with suitable choices of the matrix $\pi$.

Now specialize to the case where $k$ is algebraically closed and $\pi_{ij}=1$ for
all $i<j$. If $n$ is even, one computes that $Z_P(R) =k$, while if $n$ is odd,
one gets $Z_P(R)= k[z^{\pm1}]$ where $z= x_1x_2^{-1} x_3x_4^{-1} \cdots
x_{n-1}^{-1} x_n$. (E.g., apply the row reduction steps given in \cite{\BrGd,
Example I.14.3(1)} to $\pi$.) Therefore, we conclude that
$$\alignedat2
\Pprim R &= \{0\}  &\qquad\quad&(n \text{\;even})  \\
\Pprim R &= \{ \langle x_1x_2^{-1} x_3x_4^{-1} \cdots x_{n-1}^{-1} x_n -\lambda
\rangle \mid \lambda\in\kx \}  &&(n \text{\;odd}).
\endalignedat$$
\enddefinition

\definition{Example 4.6} Now take $R= k[x_1,\dots,x_n]$ to be a polynomial
algebra, and again choose an antisymmetric matrix $\pi= (\pi_{ij}) \in M_n(k)$.
The Poisson bracket on $k[x_1^{\pm1},\dots,x_n^{\pm1}]$ described in (4.1)
restricts to a Poisson bracket on $R$. Also, the action of the torus $H= (\kx)^n$
given by (4.2) restricts to a rational action on $R$ by Poisson automorphisms.
It is easily checked that the prime $H$-ideals in $R$ are the ideals
$$J(X)= \langle x_i\mid i\in X\rangle \qquad\quad \text{for\ } X\subseteq
\{1,\dots,n\},$$
and they are all Poisson ideals. Thus, $R$ has precisely $2^n$ prime Poisson
$H$-ideals. Theorem 4.3 therefore implies that $R$ satisfies the Poisson
Dixmier-Moeglin equivalence, and Theorem 4.1 implies that $\Pprim R$ is a
topological quotient of $\max R$. The latter result is implicit in \cite{\OPS}.

The localizations $R_J$ of Theorem 4.2 can be identified with Poisson subalgebras
of $k[x_1^{\pm1},\dots,x_n^{\pm1}]$ as follows:
$$R_{J(X)} = k[x_i^{\pm1} \mid i\in \{1,\dots,n\} \setminus X].$$
Theorems 4.2 and 4.3 imply that the Poisson primitive ideals of $R$ can be
obtained by pulling back the ideals extended from the maximal ideals of the
Poisson centers $Z_P(R_{J(X)})$. These Poisson centers can be computed as in
Example 4.5.

Now specialize to the case where $k$ is algebraically closed, $n=3$, and
$\pi_{ij}=1$ for all $i<j$. Taking account of Example 4.5, we find that the
Poisson primitive ideals of $R= k[x_1,x_2,x_3]$ in this case can be listed as
follows:
$$\alignat3
 &\langle x_1,x_2,x_3 \rangle  &\qquad\qquad\qquad\qquad&\langle
x_1-\alpha,x_2,x_3 \rangle &\qquad&(\alpha\in\kx)  \\
 &\quad\langle x_1\rangle  &&\langle x_1,x_2-\beta,x_3 \rangle &&(\beta\in\kx) 
\\
 &\quad\langle x_2\rangle  &&\langle x_1,x_2,x_3-\gamma \rangle &&(\gamma\in\kx) 
\\
 &\quad\langle x_3\rangle  &&\langle x_1x_3-\lambda x_2\rangle &&(\lambda\in\kx).
\endalignat$$
\enddefinition

\definition{Example 4.7} On the polynomial algebra $R= k[a,b,c,d]$, there is a
Poisson bracket arising from what is called the ``semiclassical limit of quantum
$2\times2$ matrices'', meaning that the standard quantized coordinate ring of
$2\times2$ matrices is a ``quantization'' of $R$ (see \cite{\BrGd, \S III.5.4},
for instance, for a discussion of this quantization concept). This Poisson bracket
is given by the following data:
$$\xalignat2
\{a,b\} &= ab  &\{b,d\} &= bd  \\
\{a,c\} &= ac  &\{c,d\} &= cd  \\
\{a,d\} &= 2bc  &\{b,c\} &= 0
\endxalignat$$
(e.g., \cite{\Vandef, Example 2.1}, \cite{\ChOh, p\. 255}, or \cite{\Oh, \S2.9},
where in the latter two papers the bracket has been scaled by 2). The torus $H=
(\kx)^4$ acts on $R$ by $k$-algebra automorphisms such that
$$\xalignat2
(\alpha_1,\alpha_2,\beta_1,\beta_2).a &= \alpha_1\beta_1a 
&(\alpha_1,\alpha_2,\beta_1,\beta_2).b &= \alpha_1\beta_2b  \\
(\alpha_1,\alpha_2,\beta_1,\beta_2).c &= \alpha_2\beta_1c 
&(\alpha_1,\alpha_2,\beta_1,\beta_2).d &= \alpha_2\beta_2d
\endxalignat$$
for $(\alpha_1,\alpha_2,\beta_1,\beta_2) \in H$. Observe that this is a
rational action by Poisson automorphisms. 

It is known that $R$ has precisely 14 prime Poisson $H$-ideals:
$$\gather
\langle a,b,c,d\rangle  \\
\langle a,b,d\rangle \qquad\quad \langle a,b,c\rangle \qquad\quad  \langle
b,c,d\rangle \qquad\quad \langle a,c,d\rangle  \\
\langle a,b\rangle \qquad\quad \langle b,d\rangle \qquad\quad \langle b,c\rangle
\qquad\quad \langle a,c\rangle \qquad\quad \langle c,d\rangle  \\
\langle b\rangle \qquad\quad \langle ad-bc\rangle \qquad\quad \langle c\rangle  \\
\langle 0\rangle.
\endgather$$
(For instance, this can be calculated as in \cite{\ChOh, pp\. 255-258}. When $k$
is algebraically closed, it also follows from \cite{\ChOh, Theorem 9} or
\cite{\Oh, p\. 2179}.) Consequently, Theorem 4.3 implies that $R$ satisfies the
Poisson Dixmier-Moeglin equivalence, recovering a result of Oh \cite{\Oh,
Proposition 2.13}. 

If $k$ is algebraically closed, one can use the information above, together with
Theorem 4.4 and calculations of appropriate Poisson centers of localizations, to
find all the Poisson primitive ideals of $R$. The calculations are essentially
the same as those employed by Cho and Oh to the same end (see \cite{\ChOh, Theorem
9} and
\cite{\Oh, p\. 2179, display before Theorem 3.4}).\quad$\diamondsuit$
\enddefinition

We conclude by sketching some examples in which Poisson primitive ideals
correspond to symplectic leaves. Although the concept of a symplectic leaf arises
in differential geometry, we shall immediately restrict attention to situations in
which it can be described in terms of algebraic geometry.

First, a {\it Poisson manifold\/} is a smooth complex manifold $M$ together with a
Poisson bracket on the algebra $C^\infty(M)$ of smooth functions on $M$. The
derivations
$\{a,-\}$ on $C^\infty(M)$ define {\it Hamiltonian vector fields\/} on $M$, and
smooth paths in $M$ whose tangent vectors come from Hamiltonian vector fields are
called {\it Hamiltonian paths\/}. One can then define the symplectic leaves of
$M$ to be the connected components relative to the relation ``connected by a
piecewise Hamiltonian path''. (Differential geometers typically prefer an
equivalent definition in terms of symplectic submanifolds -- e.g., see
\cite{\CaWe, \S5.1} or \cite{\BrGd, \S III.5.2}.) Now any smooth complex affine
variety $V$ has a natural smooth manifold structure, and a Poisson bracket on the
coordinate ring
$\O(V)$ extends uniquely to $C^\infty(V)$, thus making $V$ into a Poisson
manifold. However, the Hamiltonian
paths in $V$ need not be algebraic curves, since they are typically defined by
exponential functions, and so the symplectic leaves in $V$ are not necessarily
algebro-geometric objects. However, when these symplectic leaves are
algebraic, in the sense that they are locally closed subvarieties, they
correspond to the Poisson primitive ideals of $\O(V)$, by a result of Brown and
Gordon \cite{\BrGr, Proposition 3.6(2)}. Namely, if $\mfrak_x$ denotes the maximal
ideal of $\O(V)$ corresponding to a point $x\in V$, then the symplectic leaf
containing $x$ coincides with the set
$$\C(x) :=\{y\in V \mid \Pcore(\mfrak_y)= \Pcore(\mfrak_x) \},  \tag4.3$$
which is called the {\it symplectic core of $\mfrak_x$\/} in \cite{\BrGr, \S3.3}.
Moreover, under these hypotheses, half of the Poisson Dixmier-Moeglin equivalence
already holds, in that the Poisson primitive ideals in $\O(V)$ are precisely the
locally closed points of $\Pspec \O(V)$ \cite{\BrGr, Proposition 3.6(2)}.

Let us define an {\it affine Poisson variety\/} (over $\CC$) to be a smooth
affine complex variety $V$ together with a Poisson bracket on the coordinate ring
$\O(V)$. (In order to define Poisson structures on projective varieties, Poisson
brackets need to be rewritten in terms of bivector fields; for instance, see
\cite{\EtSc, \S1.4.1} or \cite{\KoSo, Section 1.1}.) The following proposition,
based on the ideas of Brown and Gordon \cite{\BrGr}, describes a geometric
setting in which our main results apply. For any subset $X\subseteq V$, let
$I(X)$ denote the defining ideal of the closure $\overline X$, that is, $I(X)$ is
the set of those functions in $\O(V)$ which vanish on $X$ (equivalently, on
$\overline X$).

\proclaim{Proposition 4.8} Let $V$ be an affine Poisson variety over $\CC$, and
$H$ an algebraic group acting morphically on $V$ via automorphisms of Poisson
varieties. There is an induced action of $H$ on $\O(V)$ by Poisson automorphisms.
Assume that $V$ has only finitely many $H$-orbits of symplectic leaves.

{\rm (a)} There are only finitely many prime Poisson $H$-ideals in 
$\O(V)$.

{\rm (b)} Now assume that all the symplectic leaves in $V$ are locally closed
subvarieties, and that $H= (\CC^\times)^r$ is a complex algebraic torus. Then the
rule
$$\bigl(\; H\text{-orbit\ } \L \text{\ of symplectic leaves}\; \bigr) \longmapsto
I\bigl(\; {\tsize\bigcup}\, \{L\in\L\}\; \bigr)$$
defines a bijection between the set of $H$-orbits of symplectic leaves in $V$ and
the set of prime Poisson $H$-ideals in $\O(V)$.
\endproclaim

\demo{Proof} Recall the symplectic cores $\C(x)$ defined in (4.3). There is a
bijection between the set of symplectic cores in $V$ and the set
$\Pprim \O(V)$, given by the rule $\C(x) \mapsto \Pcore(\mfrak_x)$, and this
induces a bijection between the set of $H$-orbits of symplectic cores in
$V$ and the set of $H$-orbits in $\Pprim \O(V)$.

(a) Since each symplectic core is a union of symplectic leaves \cite{\BrGr,
Proposition 3.6(1)}, it follows from our hypotheses that there are only finitely
many $H$-orbits of symplectic cores in $V$, and thus only finitely many
$H$-orbits in $\Pprim
\O(V)$. Now for any $P\in
\Pprim \O(V)$, the ideal $(P:H)$, which we shall call the {\it $H$-core\/} of $P$,
equals the intersection of the $H$-orbit $\{h(P) \mid h\in H\}$. Hence, there are
only finitely many $H$-cores of Poisson primitive ideals in $\O(V)$.

Let $Q$ be any prime Poisson $H$-ideal in $\O(V)$. By Lemma 1.1(e), $Q$ is an
intersection of Poisson primitive ideals $P_\alpha$, and since $Q$ is stable under
$H$, it is also the intersection of the $H$-cores of the $P_\alpha$.
Since there are only finitely many $H$-cores of Poisson primitive ideals in
$\O(V)$, there can only be finitely many intersections of such ideals, and
therefore part (a) is proved. 

(b) As noted above, under the present hypotheses, the symplectic leaves in $V$
are precisely the sets $\C(x)$ \cite{\BrGr, Proposition 3.6(2)}. Moreover, by
\cite{\BrGr, Lemma 3.5}, $\Pcore(\mfrak_x)= I(\C(x))$ for all $x\in V$. Due to
part (a) and the assumption that $H$ is a torus, Theorem 4.4(b) shows that the
$H$-orbits in $\Pprim \O(V)$ coincide with the $H$-strata. 

Let $\Hsymp V$ denote the set of $H$-orbits of symplectic leaves in $V$, and set
$$\theta(\L) := I\bigl(\; {\tsize\bigcup}\, \{L\in\L\}\; \bigr)$$
for $\L\in \Hsymp V$.
Since $\L$ is the $H$-orbit of a symplectic leaf $\C(x)$, for some $x\in V$, we
have
$$\theta(\L)= \bigcap_{h\in H} I\bigl( \C(h.x) \bigr)= \bigcap_{h\in H}
\Pcore(\mfrak_{h.x})= (\Pcore(\mfrak_x):H),  \tag4.4$$
and thus, by Lemma 3.1, $\theta(\L)$ is a prime Poisson $H$-ideal in $\O(V)$.
To prove (b), it remains to establish the following:
\roster
\item"(*)" For each prime Poisson $H$-ideal $Q$ in $\O(V)$, there is a unique
$\L\in \Hsymp V$ such that $\theta(\L) = Q$. 
\endroster

As in the proof of (a), there are only finitely many $H$-cores of Poisson
primitive ideals in $\O(V)$, say $(P_1:H),\dots,(P_m:H)$, for some
$P_1,\dots,P_m$ in $\Pprim \O(V)$. As also noted there, any prime
Poisson $H$-ideal $Q$ in $\O(V)$ is an intersection of $H$-cores of Poisson
primitive ideals, say $Q= (P_{i_1}:H)\cap \cdots\cap (P_{i_t}:H)$. Since $Q$ is
prime, it follows that $Q= (P_j:H)$ for some $j\in \{i_1,\dots,i_t\}$. Now $P_j=
\Pcore(\mfrak_x)$ for some $x\in V$, and thus $Q= \theta(\L)$ by (4.4), where
$\L$ is the $H$-orbit of $\C(x)$. 

Finally, suppose that also $Q= \theta(\L')= (\Pcore(\mfrak_y):H)$, where $\L'$ is
the $H$-orbit of some symplectic leaf $\C(y)$. Then $\Pcore(\mfrak_x)$ and
$\Pcore(\mfrak_y)$ lie in the same $H$-stratum $\Pprim_Q \O(V)$, and thus in the
same $H$-orbit. As noted at the beginning of the proof, this implies that $\C(x)$
and $\C(y)$ lie in the same $H$-orbit of symplectic leaves of $V$, that is,
$\L'=\L$. Therefore $\L$ is unique, and (*) is verified. \qed\enddemo

\definition{Example 4.9} Let $G$ be a connected semisimple complex Lie group,
with opposite Borel subgroups $B^{\pm}$ and corresponding Cartan subgroup $H=
B^+\cap B^-$, and let $H$ act on $G$ by left translation. There is a ``standard''
$H$-invariant Poisson structure on $G$ (e.g.,
\cite{\HoLeA, Section A.1} or \cite{\KoSo, Section 5.3}), and there are only
finitely many
$H$-orbits of symplectic leaves in $G$ \cite{\HoLeA, Theorem A.2.1}. Hence, if we
put the corresponding Poisson bracket on $\O(G)$, Proposition 4.8(a) tells us that
$\O(G)$ has only finitely many prime Poisson $H$-ideals. Therefore, $\O(G)$
satisfies the Poisson Dixmier-Moeglin equivalence, by Theorem 4.3. It also follows
from \cite{\HoLeA, \S\S A.1, A.2} that the symplectic leaves in $G$ are locally
closed subvarieties, and that the unions of $H$-orbits of symplectic leaves
coincide with the double Bruhat cells $B^+w_+B^+ \cap B^-w_-B^-$, for elements
$w_{\pm}$ in the Weyl group of $G$. Hence, Proposition 4.8(b) implies that the
prime Poisson $H$-ideals in $\O(G)$ are the defining ideals of the closures of the
double Bruhat cells in $G$.

As a specific example, let $G= SL_n(\CC)$, and take $B^+$, $B^-$, $H$ to be the
respective subgroups of upper triangular, lower triangular, diagonal matrices in
$G$. For $i,j=1,\dots,n$, let $X_{ij}\in
\O(G)$ denote the function that takes matrices to their $i,j$-entries. The Poisson
bracket on $\O(G)$ is determined by the folllowing data:
$$\{X_{ij},X_{lm}\} = \cases X_{ij}X_{lm} &\quad (i=l,\; j<m) \\
X_{ij}X_{lm} &\quad (i<l,\; j=m) \\
0 &\quad (i<l,\; j>m) \\
2X_{im}X_{lj} &\quad (i<l,\; j<m). \endcases  \tag4.5$$
Note that when $n=2$ and $\O(SL_2(\CC))$ is identified with $k[a,b,c,d]/ \langle
ad-bc-1\rangle$, the Poisson bracket above is induced from the one discussed in
Example 4.7. Since the Weyl group of $G$ can be identified with the subgroup of
signed permutation matrices in $G$, the double Bruhat cells in $G$ are easy to
identify.  Results of Fulton \cite{\Ful, p.~390} allow one to characterize these
cells in terms of vanishing and nonvanishing of certain minors, as in \cite{\BGY,
Example 4.4}.\quad$\diamondsuit$
\enddefinition

\definition{Example 4.10} Let $G$ be a connected reductive complex Lie group,
with opposite Borel subgroups $B^{\pm}$ and corrresponding Cartan subgroup $H=
B^+\cap B^-$, and let $H$ act on $G$ by left translation.  Let $P_J^+$ be a
standard parabolic subgroup containing $B^+$. There is a standard Poisson
structure on the flag variety $G/P_J^+$ (induced from the standard Poisson
structure on $G$)
\cite{\GoYa, Proposition 1.3}, which restricts to the open $B^-$-orbit $B^-.P_J^+
\subseteq G/P_J^+$ (that is, the orbit of the coset $P_J^+= eP_J^+$ under left
translation by $B^-$). These Poisson structures are invariant under the induced
left actions by $H$. According to \cite{\BGY, Theorem 1.9} and \cite{\GoYa,
Theorem 1.5}, the symplectic leaves in $G/P_J^+$ are locally closed subvarieties,
and there are only finitely many $H$-orbits of symplectic leaves. Both properties
are inherited by the affine Poisson variety $B^-.P_J^+$. Therefore, Proposition
4.8(a) shows that $\O(B^-.P_J^+)$ has only finitely many prime Poisson
$H$-ideals, and Theorem 4.3 implies that $\O(B^-.P_J^+)$ satisfies the Poisson
Dixmier-Moeglin equivalence. 

Some explicit examples of the above Poisson structures are given in \cite{\GoYa,
\S\S5.4--5.7}, and the main object of \cite{\BGY} is of this type. For the
latter, choose positive integers $m$ and $n$, let $G= GL_{m+n}(\CC)$, and take
$$P_J^+= \left\{ \left[ \matrix a&b\\ 0&c \endmatrix \right] \biggm| a\in
GL_n(\CC),\; b\in M_{n,m}(\CC),\; c\in GL_m(\CC) \right\}.$$
(As in Example 4.9, we take $B^+$, $B^-$, $H$ to be the
respective subgroups of upper triangular, lower triangular, diagonal matrices in
$G$.) In the present case, $B^-.P_J^+$ is isomorphic, as a
Poisson variety, to the matrix variety $M_{m,n}(\CC)$, equipped with its standard
Poisson structure \cite{\BGY, Proposition 3.4}. (See \cite{\BGY, \S1.5} for the
standard Poisson structure on $M_{m,n}(\CC)$. If we take $X_{ij}\in
\O(M_{m,n}(\CC))$, for $i=1,\dots,m$ and $j=1,\dots,n$, to be the usual
matrix-entry functions, the data for the Poisson bracket on $\O(M_{m,n}(\CC))$
are given by (4.5).) The $H$-orbits of symplectic leaves in $M_{m,n}(\CC)$ are
described in three different ways in \cite{\BGY, Theorems 3.9, 4.2, 5.11}.
Combining this information with Proposition 4.8(b) yields descriptions of the
prime Poisson $H$-ideals in $\O(M_{m,n}(\CC))$.\quad$\diamondsuit$
\enddefinition


\Refs

\widestnumber\key{{\bf 99}}

\ref\no\BrGd \by K. A. Brown and K. R. Goodearl \book Lectures on Algebraic
Quantum Groups \bookinfo Advanced Courses in Math. CRM Barcelona \publaddr Basel
\yr 2002 \publ Birkh\"auser \endref

\ref\no\BGY \by K. A. Brown, K. R. Goodearl, and M. Yakimov \paper Poisson
structures on affine spaces and flag varieties. I. Matrix affine Poisson space
\jour Advances in Math. \finalinfo in press, posted at
arxiv.org/abs/math.QA/{\allowlinebreak}0501109 \endref 

\ref\no\BrGr \by K. A. Brown and I. Gordon \paper Poisson orders, symplectic
reflection algebras and representation theory \jour J. reine angew. Math. \vol
559 \yr 2003 \pages 193-216 \endref

\ref\no\CaWe \by A. Cannas da Silva and A. Weinstein \book Geometric Models for
Noncommutative Algebras \bookinfo Berkeley Math. Lecture Notes 10 \publ Amer.
Math. Soc. \yr 1999 \publaddr Providence \endref

\ref\no\ChOh \by E. H. Cho and S.-Q. Oh \paper Primitive spectrum of quantum
$2\times2$ matrices and associated Poisson structure \jour Far East J. Math.
Sci. \vol 6 \yr 1998 \pages 251-259 \endref

\ref\no\Dix \by J. Dixmier \paper Id\'eaux primitifs dans les alg\`ebres
enveloppantes \jour J. Algebra \vol 48 \yr 1977 \pages 96--112 \endref

\ref\no \Dixbook \bysame \book Enveloping Algebras \publ North-Holland \publaddr
Amsterdam
\yr 1977 \endref

\ref\no\EtSc \by P. Etingof and O. Schiffmann \book Lectures on Quantum Groups
\publaddr Boston \yr 1998 \publ Internat. Press \endref

\ref\no\Far \by D. R. Farkas \paper Characterizations of Poisson algebras \jour
Communic. in Algebra \vol 23 \yr 1995 \pages 4669-4686 \endref

\ref\no\FaLe \by D. R. Farkas and G. Letzter \paper Ring theory from symplectic
geometry \jour J. Pure Applied Algebra \vol 125 \yr 1998 \pages 155-190 \endref

\ref\no\Ful \by W. Fulton \paper Flags, Schubert polynomials,
degeneracy loci, and determinantal formulas \jour Duke Math. J.
\vol 65 \yr 1991 \pages 381-420 \endref

\ref\no \specstrat \by K. R. Goodearl and E. S. Letzter\paper The Dixmier-Moeglin
equivalence in quantum coordinate rings and quantized Weyl algebras \jour Trans.
Amer. Math. Soc. \vol 352 \yr 2000 \pages 1381-1403 \endref

\ref\no\qaffquo \bysame \paper Quantum $n$-space as a quotient of classical
$n$-space \jour Trans. Amer. Math. Soc. \vol 352 \yr 2000 \pages 5855-5876 \endref

\ref\no\GoYa \by K. R. Goodearl and M. Yakimov \paper Poisson structures on
affine spaces and flag varieties. II \paperinfo preprint 2005 \finalinfo posted at
arxiv.org/abs/math.QA/0509075
\endref


\ref\no\HoLeA \by T. J. Hodges and T. Levasseur\paper Primitive ideals of
${\bold C}_q[SL(3)]$ \jour Commun. Math. Phys. \vol 156 \yr 1993 \pages 581-605
\endref

\ref\no\HoLeB \bysame \paper Primitive ideals of ${\bold
C}_q[SL(n)]$ \jour J. Algebra \vol 168 \yr 1994 \pages 455-468 \endref

\ref\no\Hum \by J. E. Humphreys \book Linear Algebraic Groups \bookinfo Corr.
Third Printing \publaddr New York \yr 1981 \publ Springer-Verlag \endref

\ref\no\Jos \by A. Joseph \book Quantum Groups and Their Primitive Ideals
\bookinfo Ergebnisse der Math. (3) 29 \publ Spring\-er-Ver\-lag \publaddr Berlin
\yr 1995
\endref

\ref\no\KoSo \by L. I. Korogodski and Ya. S. Soibelman \book Algebras of Functions
on Quantum Groups: Part I \bookinfo Math. Surveys and Monographs 56 \publ Amer.
Math. Soc. \publaddr Providence \yr 1998
\endref 

\ref\no\Moe \by C. Moeglin \paper Id\'eaux primitifs des alg\`ebres
enveloppantes \jour J. Math. Pures Appl. \vol 59 \yr 1980 \pages
265--336 \endref

\ref\no \Oh \by S.-Q. Oh \paper Symplectic ideals of Poisson algebras and the
Poisson structure associated to quantum matrices \jour Comm. Algebra \vol 27 \yr
1999 \pages 2163-2180  \endref

\ref\no\OhPa \by S.-Q. Oh and C.-G. Park \paper Primitive ideals in the coordinate
ring of quantum Euclidean space \jour Bull. Austral. Math. Soc. \vol 58 \yr
1998 \pages 57-73 \endref

\ref\no\OPS \by S.-Q. Oh, C.-G. Park, and Y.-Y. Shin \paper Quantum $n$-space and
Poisson $n$-space \jour Comm. Alg. \vol 30 \yr 2002 \pages 4197-4209 \endref

\ref\no\Vandef \by M. Vancliff \paper The defining relations of quantum $n\times
n$ matrices\jour J. London Math. Soc. (2) \vol 52\yr 1995\pages 255-262 \endref

\ref\no\Van \bysame
\paper Primitive and Poisson spectra of twists of
polynomial rings \jour Algebras and Representation Theory \vol 2 \yr 1999 \pages
269-285 \endref

\endRefs
\enddocument